\documentclass[12pt]{amsart}
\usepackage{graphicx,amsfonts}
\usepackage{amsmath,amssymb}
\usepackage[all]{xy}

\newcommand{\C}{\mathbb{C}}

\newcommand{\Z}{\mathbb{Z}}
\newcommand{\Q}{\mathbb{Q}}
\newcommand{\A}{\mathcal{A}}
\newcommand{\K}{\mathcal{K}}
\newcommand{\z}{\zeta}

\newcommand{\id}{\mathop{\rm id}\nolimits}
\renewcommand{\phi}{\varphi}
\renewcommand{\k}{\varkappa}

\newcommand\Am{{\mathcal A}^h(m)}
\newcommand\Amh{\hat{\mathcal A}^h(m)}
\newcommand\At{\hat{\mathcal A}^h(3)}

\newcommand{\ig}{\includegraphics}
\newcommand{\rb}{\raisebox}

\newtheorem*{conjecture}{Conjecture}
\newtheorem*{theorem}{Theorem}

\title{Conway polynomial and Magnus expansion}
\thanks{Supported by grants RFBR 08-01-00379-a, NSh 709.2008.1, 
JSPS S-09018.}
\author{S.\,V.\,Duzhin}

\begin{document}
\begin{abstract}
The Magnus expansion is a universal finite type invariant of pure braids
with values in the space of horizontal chord diagrams. The Conway polynomial
composed with the short-circuit map from braids to knots gives rise to a
series of finite type invariants of pure braids and thus factors through
the Magnus map. We describe explicitly the resulting mapping from horizontal
chord diagrams on 3 strands to univariate polynomials and evaluate it on
the Drinfeld associator obtaining, conjecturally, a beautiful generating 
function whose coefficients are alternating sums of multiple zeta values.
\end{abstract}

\maketitle

\section{Introduction}
\label{sec-intr}

We assume that the reader is familiar with fundamentals of knot theory,
braid groups and finite type (Vassiliev) invariants. All these preliminaries
can be found, for instance, in \cite{PS}.

The short-circuit closure of pure braids \cite{MSt} induces a map from pure
braids onto the set of (topological types of) oriented knots.  Any Vassiliev
invariant of knots thus becomes a finite-type invariant of pure braids. 
There is a universal finite type invariant of pure braids given by the
Magnus expansion.  In this paper, we will explicitly describe the map from
horizontal chord diagrams on 3 strands obtained by factoring the Conway
polynomial pulled back to pure 3-braids through Magnus expansion.  The
result is described by a peculiar combinatorial map from ordered partitions
of an integer into non-ordered partitions of the same integer.

In Section \ref{sec-pbme} we say introductory words about the group of pure
braids and give the construction of the Magnus expansion.  Section
\ref{sec-scc} is devoted to the construction of the short-circuit closure
relating braids to knots. In Section \ref{sec-scp} we speak about the
Conway polynomial of braids transferred from knots via short-circuit closure
and state the main theorem, whose proof is given in Section \ref{sec-proof}.
In Section \ref{sec-assoc} we try to evaluate the mapping obtained in the
main theorem, on
the Drinfeld associator and state the results of our computer calculations
and the corresponding conjecture. Finally, Section \ref{sec-prob} lists some
open problems related to the material of the paper.

I am indebted to Jacob Mostovoy who read the first version of the paper
and made numerous useful remarks.

\section{Pure braids and Magnus expansion}
\label{sec-pbme}

Let $P_m$ be the group of pure braids on $m$ enumerated vertical strands
with multiplication defined as the concatenation from top to bottom.
It is generated by the elements $x_{ij}$, $1\leq i<j\leq m$,
representing one full positive twist between the $i$-th and the $j$-th
strands with all the remaining strands of the braid being strictly
vertical and placed behind these two:
$$
x_{ij}=\rb{-12mm}{\ig[height=25mm]{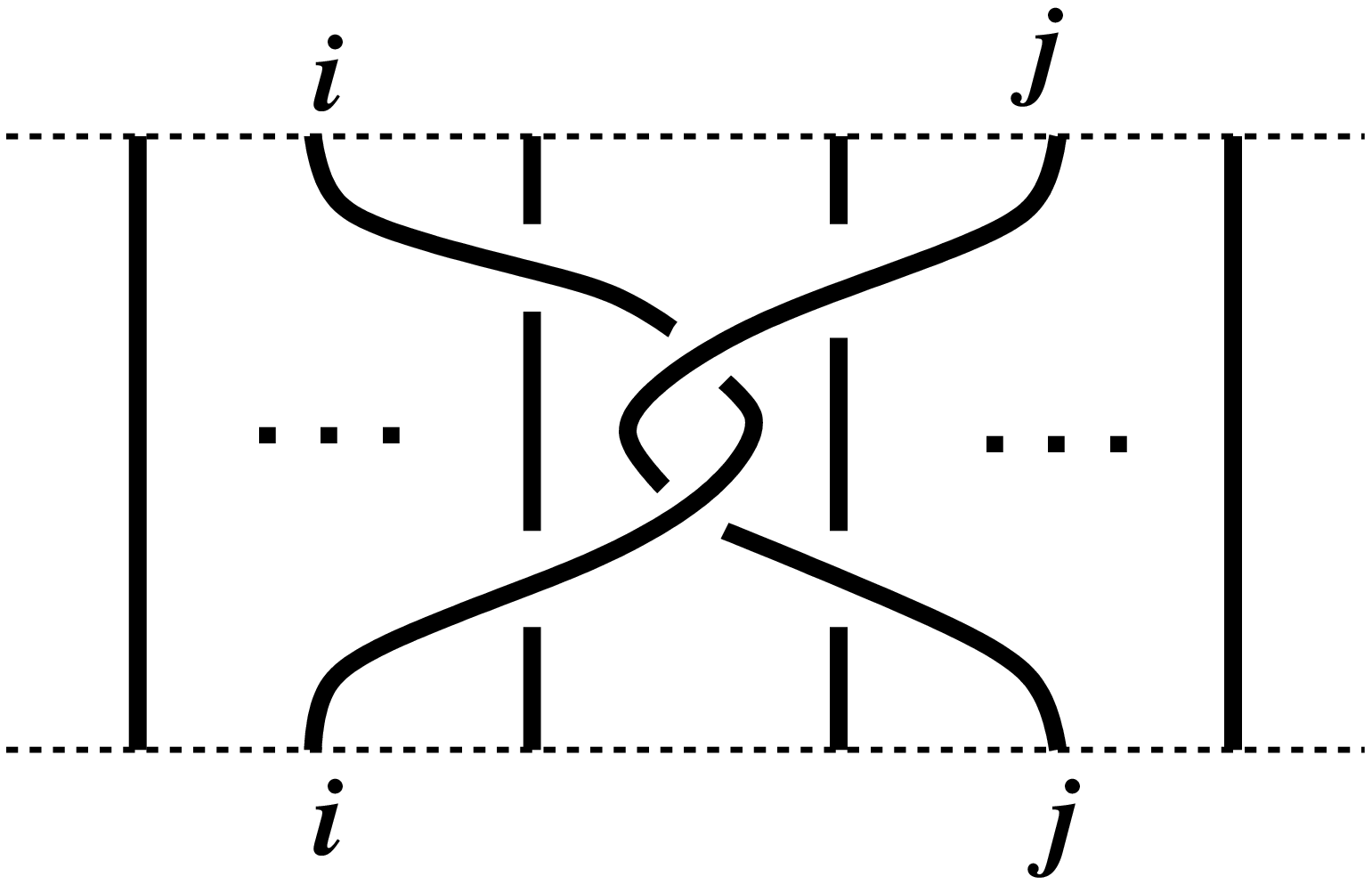}}
$$
Defining relation between these generators can be found, for instance,
in \cite{Bir,Dr}; we will not need them here.

There is a semi-direct decomposition \cite{Bir}
$$
  P_{m}\cong F_{m-1}\ltimes\ldots F_{2}\ltimes F_1,
$$
where $F_k$ is a free group on $k$ generators, implemented in our case
as the subgroup of $P_m$ generated by the set
$x_{1,k+1},x_{2,k+1},\dots,x_{k,k+1}$.
This decomposition assures that any pure braid
can be uniquely written in the \textit{combed form}
$\prod_s x_{i_sj_s}^{a_s}$ with $j_1\geq j_2\geq\dots$.
where $a_s$ are arbitrary nonzero integers and in this product no two identical 
generators follow each other (that is, the word is reduced).

The \textit{Magnus expansion} is a map from 
$P_m$ into the $\Z$-algebra of formal power series in $\binom{m}{2}$
non-commuting variables $t_{ij}$, $1\leq i<j\leq m$, defined by 
$$
\mu_m(\beta) =  \prod_s (1+t_{i_sj_s})^{a_s},
$$
if $\prod_s x_{i_sj_s}^{a_s}$ is the combed form of the braid $\beta$.
Here the negative powers are understood as usual, according to the rule
$(1+t)^{-1}=1-t+t^2-t^3+\dots$, --- this is why we need power series,
not just polynomials, in the construction of $\mu_m$.
\smallskip

{\bf Example.} To compute the value $\mu_3(x_{12}x_{23})$, we first find the 
combed form of this braid 
$$
 x_{12}x_{23} = x_{13}x_{23}x_{13}^{-1}x_{12}
$$
and then write:
\begin{eqnarray*}
\mu_3(x_{12}x_{23}) &=&
(1+t_{13})(1+t_{23})(1-t_{13}+t_{13}^2-\dots)(1+t_{12})\\
&=& 1 + t_{12}+t_{23} + t_{13}t_{23} -t_{23}t_{13} +t_{23}t_{12} +...
\end{eqnarray*}

From a broader perspective, it makes sense to consider
the codomain of the mapping $\mu_m$ as the completed quotient of the
algebra $\Z[\{{t_{ij}\}}_{1\le i < j\le m}]$
over the ideal generated by the elements
$[t_{ij},t_{kl}]$ and $[t_{ij},t_{ik}+t_{jk}]$
where all indices are supposed to be distinct and
$t_{pq}$ is understood as $t_{qp}$ if $p>q$.
We denote this algebra by $\Am$ and view its generating monomials
as horizontal chord diagrams on $m$ vertical strands (each variable $t_{ij}$
represents a chord connecting the $i$-th and $j$-th strands;
the product of variables is understood as vertical concatenation from top to
bottom). Example:
$$
  t_{13}t_{23}^2t_{12} = \rb{-8mm}{\ig[height=17mm]{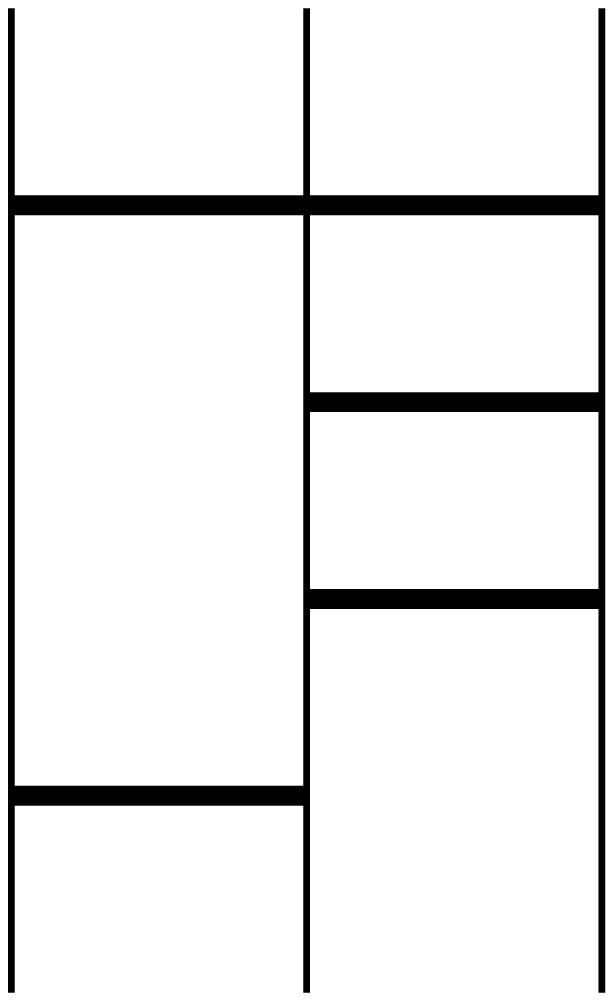}}\ .
$$
We will denote the completion of $\Am$, that is the corresponding algebra 
of formal series, by $\Amh$.

Call a horizontal chord diagram \textit{descending}, if it is represented
by a monomial $\prod_s t_{i_sj_s}^{a_s}$ satisfying $j_1\geq j_2\geq\ldots$
By definition, the set of descending diagrams is in one-to-one 
correspondence with the set of positive combed braids $P_m^+$ (braids whose 
combed form contains only positive powers of the generators $x_{ij}$).
The set of descending chord diagrams forms a basis of the free abelian group
$\Am$ (see \cite{MW}, Sec. 3-2)\footnote{In a
more general setting, this fact easily follows from Theorem 3.1 of \cite{Pap}.
It can also be proved by using non-commutative Gr\"obner
bases (I thank A.\,Khoroshkin for teaching me the idea of this proof).
Closely related formulations
are scattered in the works by T.\,Kohno, V.\,Drinfeld, An.\,Kirillov,
S.\,Yuzvinsky etc.}, 
therefore we have a module isomorphism $\Z P_m^+\cong\Am$.

By an \textit{invariant of braids}, we understand any mapping from the braid
group $P_m$ into an arbitrary set --- we are interested only in its
invariance under the braid isotopy, tacitly assumed in the definition of
$P_m$, and not in the invariance under the renumbering of strands etc.  For
pure braids, just like in the classical case of knots, one can define the
notion of finite type (Vassiliev) invariants, see \cite{BN,MW,CDM}.  It
turns out that the Magnus expansion truncated to any degree $n$ is a
Vassiliev invariant of order $n$.  Moreover, the following theorem holds:

\begin{theorem} {\rm(}\cite{MW,Pap,CDM}{\rm)}
The mapping $\mu_m:P_m\to\Amh$ is a universal finite type invariant of
pure braids in the sense that for any degree $n$ invariant $f:P_m\to\Z$
there exists a map $g:\Amh\to\Z$ vanishing on all monomials of degree 
greater than $n$ and such that $f=g\circ\mu_m$.
\end{theorem}

\textbf{Remark.}\label{magnus_gen}
In fact, one may define a universal finite type invariant
of pure braids by sending each $x_{ij}$ in the combed form into
$1+c_{ij}t_{ij}+T_{ij}$ where $c_{ij}$ are any nonzero constants and 
$T_{ij}$ are arbitrary series starting with
degree greater than one. A remarkable instance of this construction
(with values in $\Amh\otimes\C$)
is provided by the Kontsevich integral (\cite{BN,CDM}). Its
advantage over the usual Magnus expansion consists in multiplicativity;
however, the definition of the Kontsevich integral is much more involved
and its computation much more difficult; moreover, its value depends
on the placement of the endpoints of the braid. 
For example, the Kontsevich integrals of the
generating braids of the group $P_3$, where the strand endpoints are
collinear and equidistant, are infinite series with
the following terms up to degree 2:
\begin{eqnarray*}
  I(x_{12}) &=& 1-A+\frac{1}{2}A^2-\frac{i\ln2}{2\pi}[B,C]+\dots,\\
  I(x_{13}) &=& 1-C+\frac{1}{2}C^2+\frac{1}{2}[A,B]+\dots,\\
  I(x_{23}) &=& 1-B+\frac{1}{2}B^2+\frac{i\ln2}{2\pi}[C,A]+\dots,\\
\end{eqnarray*}
where $A=t_{12}, B=t_{23}, C=t_{13}$ and we remind that $[A,B]=[B,C]=[C,A]$
according to the definition of $\At$.
The reader may wish, by way of exercise, to check that these relations agree
with the commutation relations in the group $P_3$ (expressed by the fact
that the element $x_{12}x_{13}x_{23}$ is central, see \cite{Dr}).

\section{Short-circuit closure}
\label{sec-scc}

Alongside with the ordinary (Artin) closure which turns braids
into links, there is another operation of closing pure braids into 
oriented knots, called \textit{short-circuit closure}, see \cite{MSt,CDM}. 
It is defined by connecting pairwise by short arcs the upper endpoints
number $2i$ and $2i+1$ and the lower endpoints number $2i-1$ and $2i$
thus obtaining a long knot with two loose endpoints. 
Attaching an additional arc, one obtains a conventional compact knot.
Orientation is chosen so that the leftmost strand of the braid is oriented 
downwards. For example:
$$
 \rb{-6mm}{\ig[width=25mm]{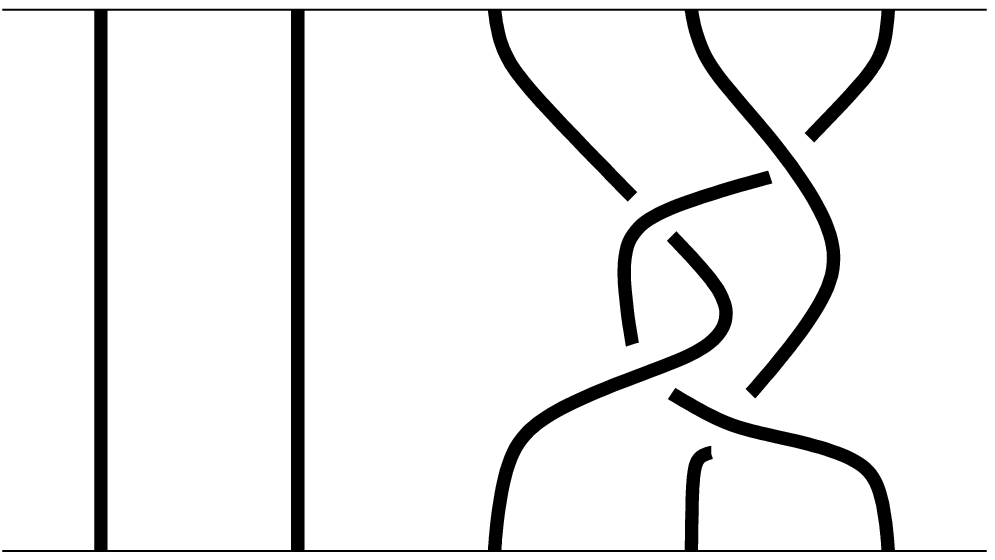}}\ \longmapsto\ 
 \rb{-10mm}{\ig[width=25mm]{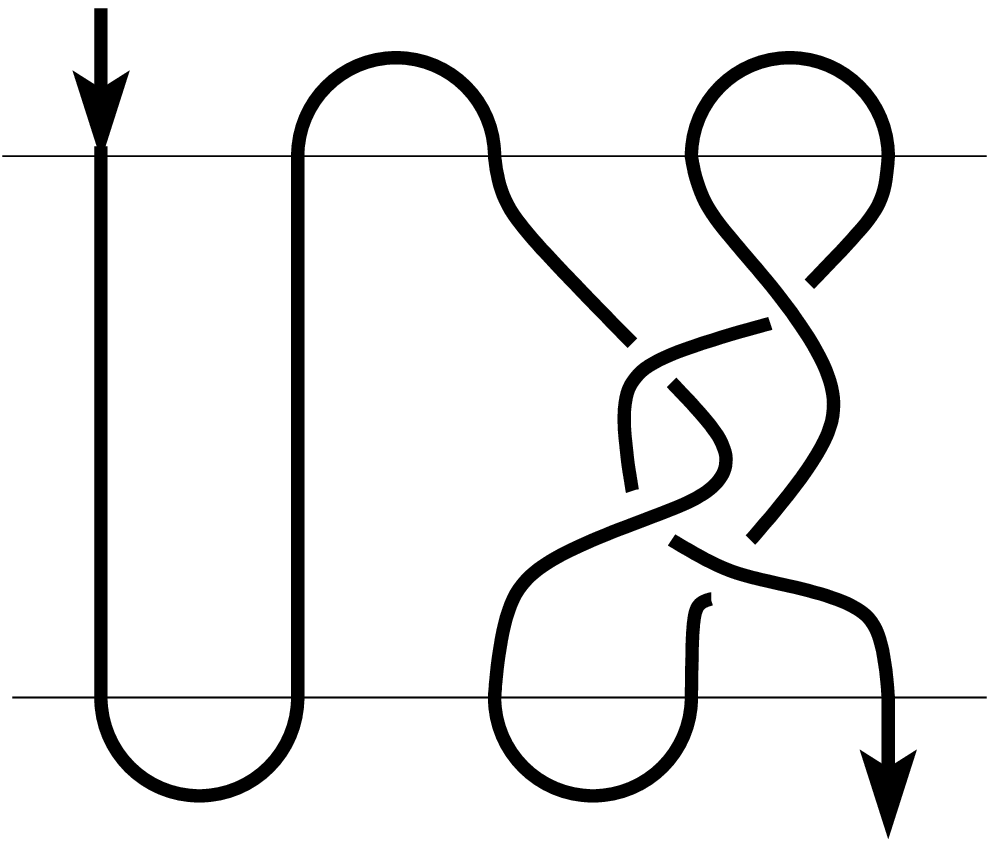}}\ \longmapsto\
 \rb{-6mm}{\ig[width=13mm]{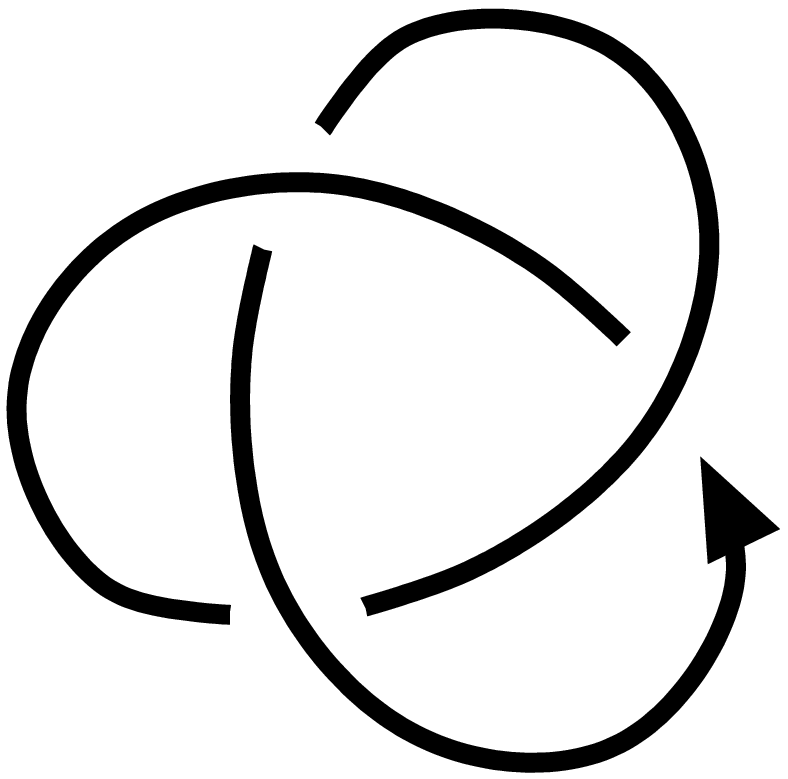}}\ .
$$
(The non-oriented version of this operation, called {\it plat closure},
was earlier introduced and studied by J.~Birman \cite{Bir}.)

It is easy to see that the short-circuit closures of braids with different
number of strands are consistent with the inclusions $P_m\to P_{m+1}$ 
(adding a vertical strand on the right), so that we obtain a well-defined 
map $\k$ from the group $P_\infty:=\cup_{m\ge1}P_m$ to the set of 
oriented knots $\K$.
A theorem of Mostovoy and Stanford asserts that this map is
onto and that it identifies $\K$ with double cosets of the group $P_\infty$
over two special subgroups, see \cite{MSt}.

In the particular case $m=3$, the image of $\k_3=\k\mid_{P_3}$
coincides with the set of all 2-bridge (rational) knots (see
\cite{Mur} for an introduction to rational knots).
Indeed, the short-circuit closure of the braid 
$x_{13}^{a_1}x_{23}^{b_1}\ldots x_{13}^{a_k}x_{23}^{b_k}$,
respectively $x_{13}^{a_1}x_{23}^{b_1}\ldots
x_{13}^{a_k}x_{23}^{b_k}x_{13}^{a_{k+1}}$,
where all $a_i$, $b_i$ are non-zero integers, is the rational 
knot corresponding to the continuous fraction with denominators
$$
(2a_1,-2b_1,...,2a_k,-2b_k+1),
$$ 
respectively, 
$$
(2a_1,-2b_1,...,2a_k,-2b_k,2a_k+1).
$$
Now, a simple number-theoretic argument shows that any rational number with an
odd denominator (those corresponding to knots, not links) has a continuous
fraction of this kind (the last number is odd, while all the previous ones
are even).

\section{Symbol of the Conway polynomial}
\label{sec-scp}

By taking composition with the short-circuit closure, any invariant 
of knots can be converted into an invariant of braids. For example, the
Conway polynomial of knots $\nabla:\K\to\Z[t]$
induces an invariant of pure braids $\nabla\circ\k_m:P_m\to\Z[t]$
(the Conway polynomial of braids).
For every $n$, the coefficient of $t^{2n}$ in this polynomial is a Vassiliev
invariant of degree $2n$.
By the universality of the Magnus expansion, there is a mapping
(the ``symbol'' of the Conway polynomial)
$\chi_m:\Amh\to\Z[t]$ such that $\nabla\circ\k_m=\chi_m\circ\mu_m$. 
We succeeded in finding an explicit description of the symbol only for pure
braids on 3 strands.

The following theorem describes the values of the
mapping $\chi=\chi_3$ on descending chord diagrams
(which form a basis of the free abelian group $\A^h(3)$). 
Let $A=t_{12}$, $B=t_{23}$, $C=t_{13}$.

\begin{theorem}
Any descending chord diagram on three strands is a (positive)
word in the letters $A,B,C$ where all $A$'s come at the end.
We claim that for any words $w$, $w_1$, $w_2$

{\rm(1)} $\chi(wA)=0$.

{\rm(2)} $\chi(Bw)=0$.

{\rm(3)} $\chi(w_1B^2w_2)=0$.

Assertions (1), (2), (3) leave us with just two kinds of words:
$$
 C^{c_1}B\cdot\ldots\cdot C^{c_{k-1}}BC^{c_k}
$$ 
and
$$
 C^{c_1}B\cdot\ldots\cdot C^{c_{k-1}}BC^{c_k}B
$$
which we encode respectively by $[c_1,\dots,c_k]$ and $[c_1,\dots,c_k]'$.

{\rm(4)} The values of $\chi$ on the elements of the second kind
are reduced to its values on the elements of the first kind:
$$
\chi([c_1,\dots,c_k]')=t^{-2}\chi([c_1,\dots,c_k,1]).
$$

It thus remains to determine $\chi$ on the elements $[c_1,...,c_k]$.

{\rm(5)} We have
$$\chi([c_1,\dots,c_k]) =
(-1)^{k-1}\bigl(\prod_{i=1}^{k-1}p_1p_{c_i-1}\bigr)\cdot p_{c_k},
$$
where $p_s=\chi([s])$ is a sequence of polynomials in $t$ that can be defined 
recursively by $p_0=1$,
$p_1=t^2$ and $p_{s+2}=t^2(p_s+p_{s+1})$ for $s\geq0$.
In particular, the value of $\chi$ on the empty chord diagram {\rm(}unit 
of the algebra $\At${\rm)} is 1.
\end{theorem}

\textbf{Remark 1.}
Note that the polynomial $p_k=\chi([k])$ is equal to $t^k\nabla(T_{k+1,2})$
where the letter $T$ denotes the torus link with given parameters (in the
case when this is a 2-component link, correct orientations of the components
must be chosen) and can be written explicitly as
$$
  p_k=\sum_{k/2\le j\le k}\binom{j}{2j-k}t^{2j}.
$$
\smallskip

\textbf{Remark 2.}
The image of $\chi$ belongs to the commutative algebra generated by
polynomials $p_1$, $p_2$ etc, whose additive basis can be identified 
with (unordered) partitions. In this setting, the map $\chi$ is defined
by a transformation of ordered partitions into unordered partitions
according to the rule
$$
  [c_1,...c_k] \mapsto (1^{k-1},c_1-1,...,c_{k-1}-1,c_k).
$$

\textbf{Examples.}
\begin{eqnarray*}
\chi(1)&=&1,\\
\chi(B)&=&0,\\
\chi(C)&=&t^2,\\
\chi(CB)&=&-t^2,\\
\chi(BC)&=&0,\\
\chi(C^3BC^3)&=&-p_1p_2p_3=-t^2(t^4+t^2)(t^6+2t^4).
\end{eqnarray*}

\section{Proof of the theorem}
\label{sec-proof}

We must find a map $\chi$ rendering commutative the diagram
$$
\xymatrix{
P_3 \ar[r]^{\mu_3} \ar@<.6ex>[d]^{\k} & 
\At \ar[d]^{\chi} \\
\K \ar[r]^{\nabla} & \Z[[t]]
}
$$

Extend the Magnus expansion linearly to the map 
$\widehat{\Z P_3}\to\At$ denoted by the same letter $\mu_3$.
We will prove the theorem by finding a left inverse of $\mu_3$, 
that is a mapping $\nu_3:\At\to\widehat{\Z P_3}$ such that $\nu_3\circ\mu_3=\id$.
Indeed, the set of decreasing chord diagrams on 3 strands is in one-to-one
correspondence with the set of positive braids $P_3^{+}$.
The correspondence is defined simply by setting
$x_{ij}\leftrightarrow t_{ij}$.
Identifying a word $w$ in $x_{ij}$ with the corresponding word in $t_{ij}$,
we see that 
$$
  \mu_3(w)=\sum_{w'\subseteq w}w'
$$
for any positive word $w$.
It is easy to check that the inverse of this map $\widehat{\Z P_3^{+}}\to\At$
is given by the formula
$$
  \nu_3(w)=\sum_{w'\subseteq w}(-1)^{|w|-|w'|}w',
$$
where the absolute value of a word denotes its length (or total exponent).
The diagram
$$
\xymatrix{
\widehat{\Z P_3} \ar[r]^{\mu_3} \ar@<.6ex>[d]^{\k} & 
\At \ar[d]^{\chi} \ar@<.6ex>[l]^{\nu_3}\\
\widehat{\Z\K} \ar[r]^{\nabla} & \Z[[t]]
}
$$
shows that we have $\chi=\nabla\circ\k\circ\nu_3$ and, consequently,
$$
 \chi(w) = \sum_{w'\subseteq w} (-1)^{|w|-|w'|} \nabla(\k(w')),
$$
where $w'$, a word in the letters $t_{ij}$, is understood,
via the mentioned identification, as a word in the letters $x_{ij}$,
that is, as a positive pure braid.

We will consecutively prove the five assertions of the theorem
by applying this expression for $\chi$ and splitting the sum over all
$2^n$ subwords $w'\subseteq w$ into appropriate subsums of 2, 4, ..., $2^k$
terms.
\smallskip

(1).
Split the sum into pairs $\pm\big(\nabla(\k(w'A))-\nabla(\k(w'))\big)$
and notice that the knots $\k(w'A)$ and $\k(w')$ are isotopic.
\smallskip

(2).
The same argument for the pairs of knots $\k(Bw')$ and $\k(w')$.
\smallskip

(3).\label{mt3}
The sum giving $\chi(w_1B^2w_2)$ consists of quadruples defined by a  
choice of subwords $w'_1\subseteq w_1$, $w'_2\subseteq w_2$:
$$
\pm\big(\nabla(\k(w'_1B^2w'_2))-2\nabla(\k(w'_1Bw'_2))
+\nabla(\k(w'_1w'_2))\big).
$$
We shall prove that any such quadruple sums to zero.

\newcommand{\kbt}{\rb{-12mm}{\ig[height=27mm]{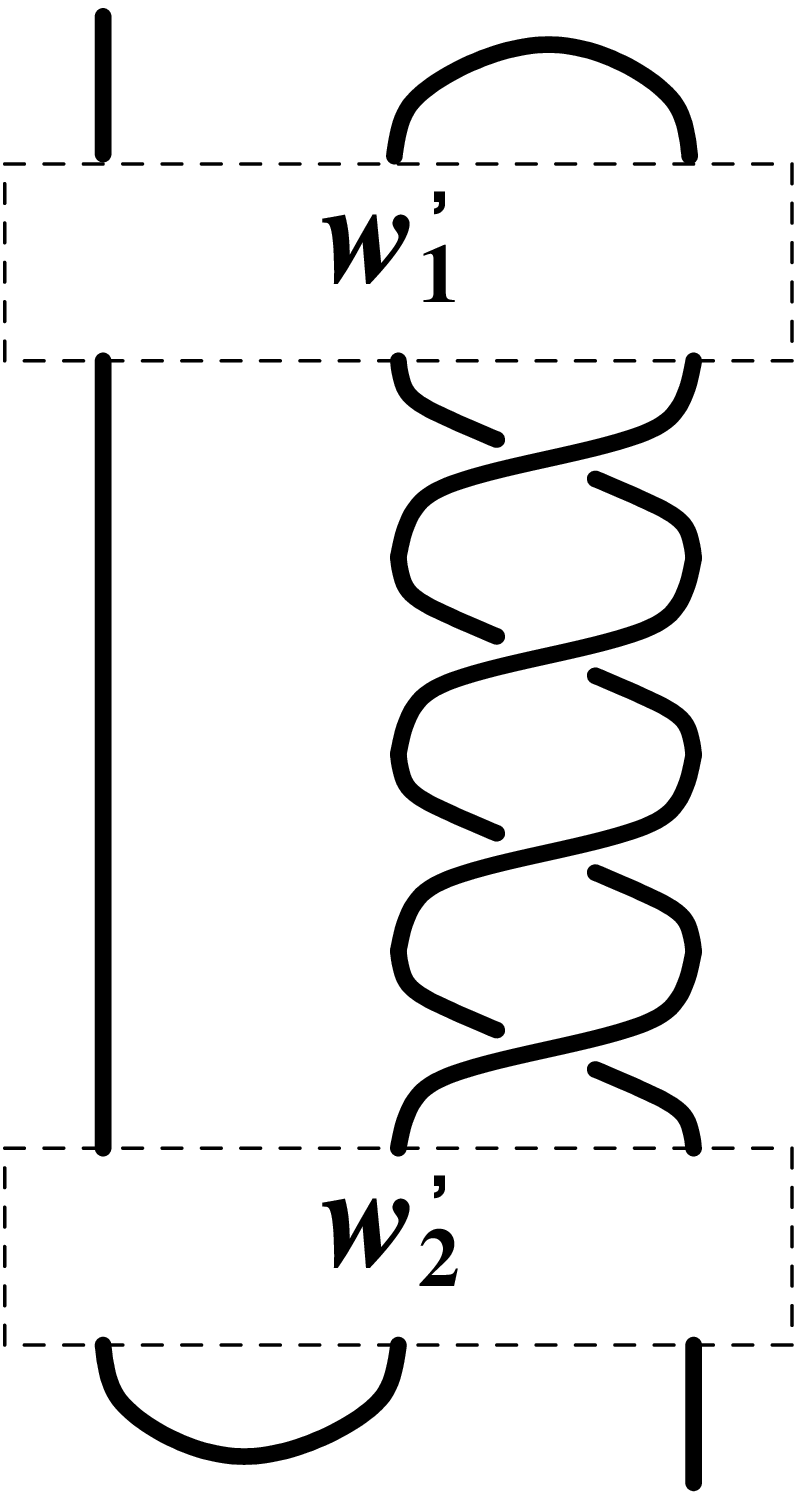}}}
\newcommand{\kbo}{\rb{-12mm}{\ig[height=27mm]{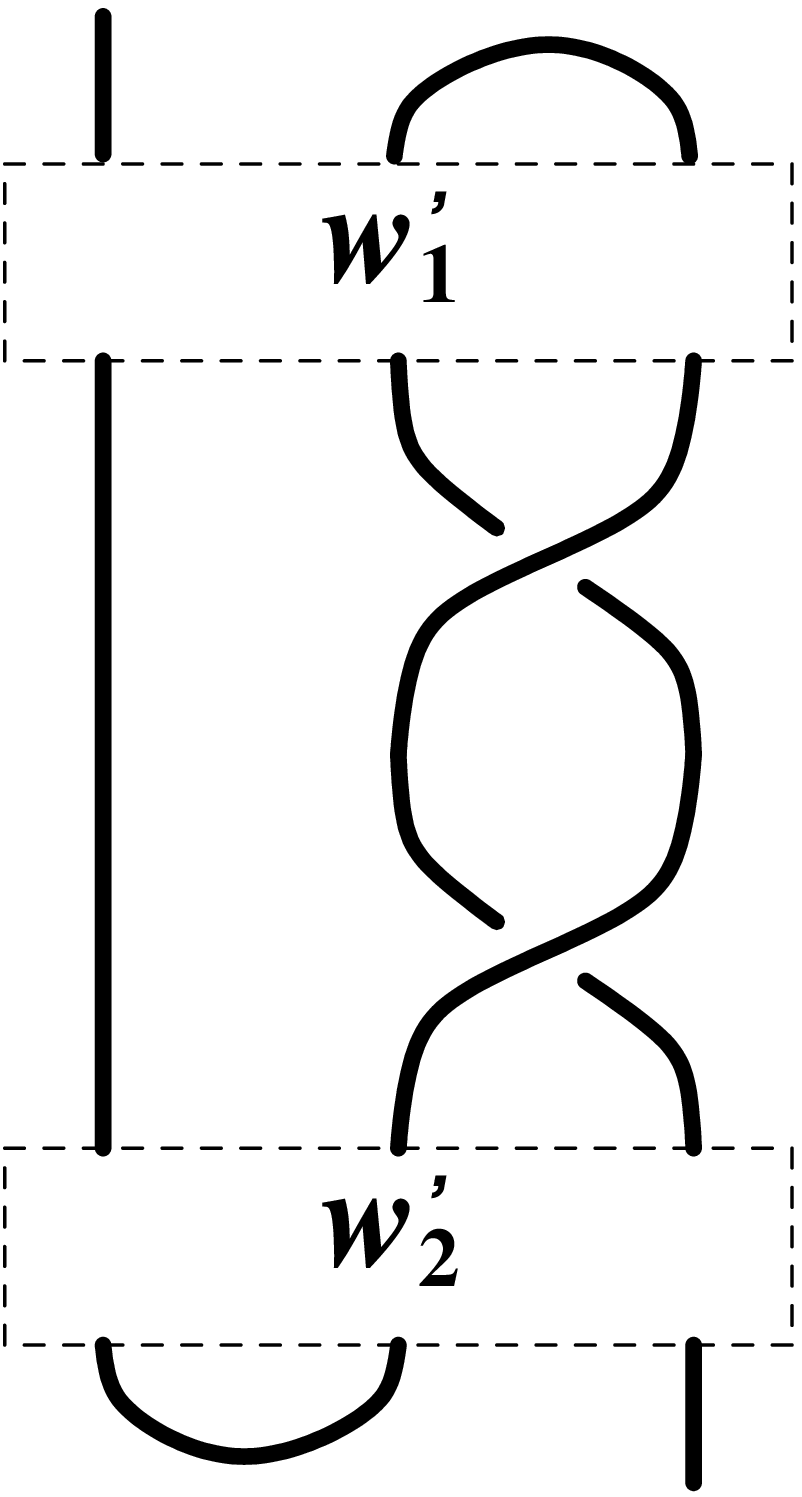}}}
\newcommand{\kbz}{\rb{-12mm}{\ig[height=27mm]{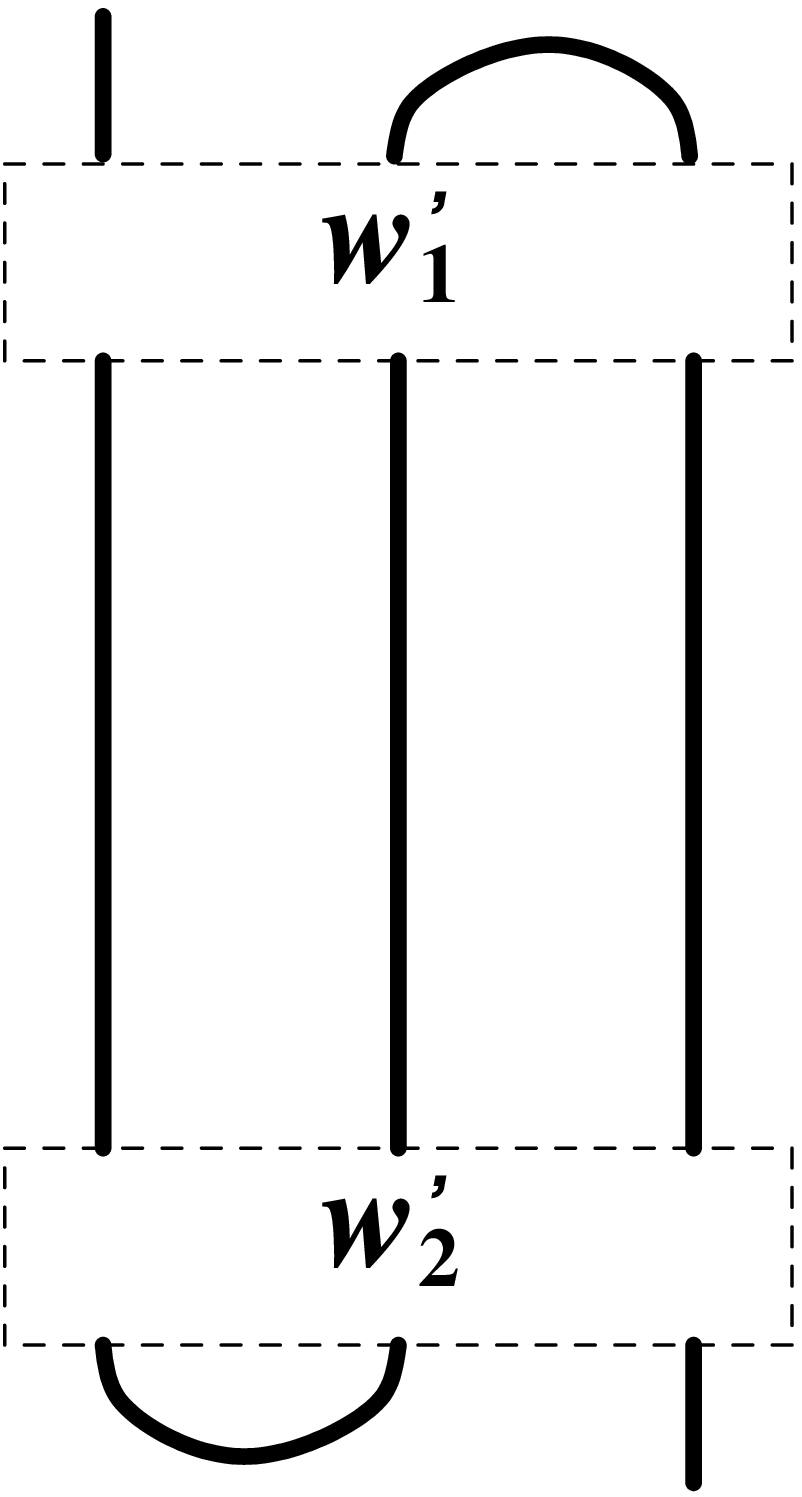}}}
\newcommand{\kbs}{\rb{-10mm}{\ig[height=21mm]{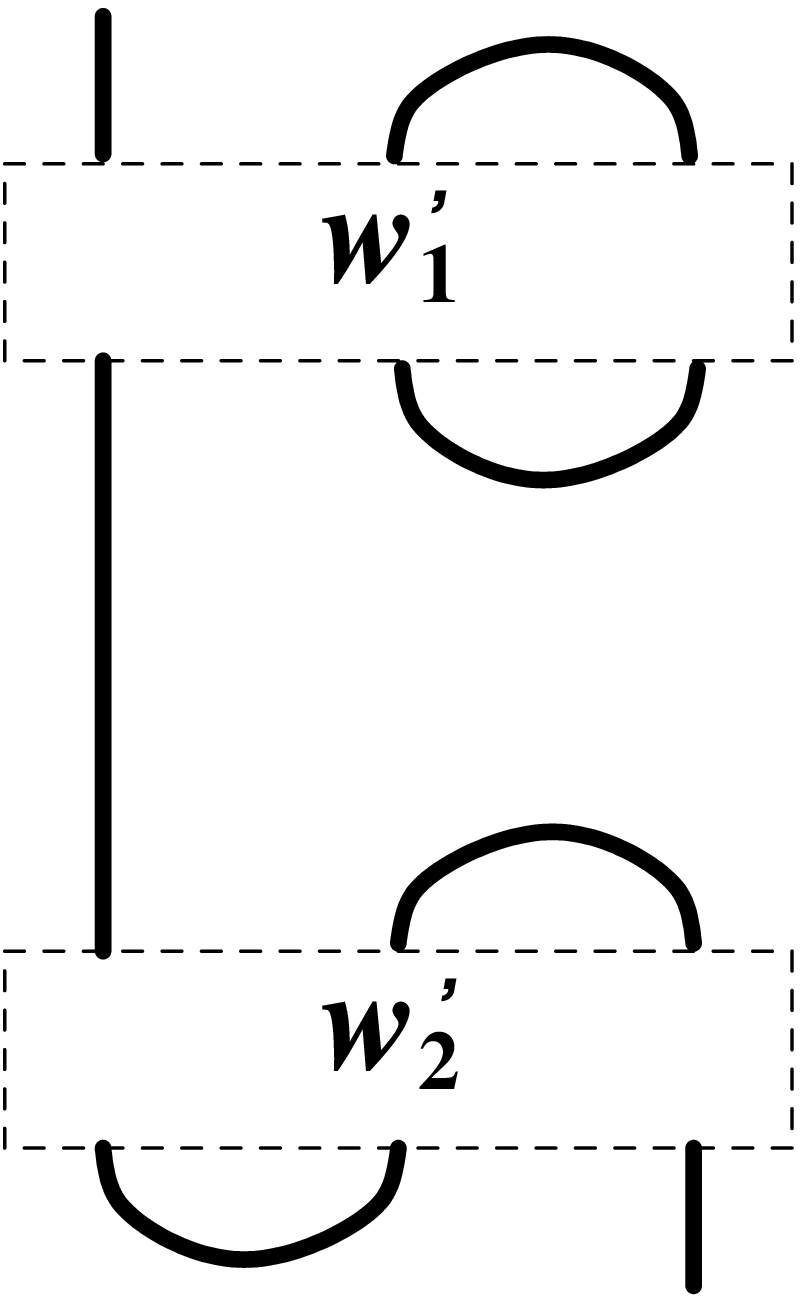}}}

Indeed, the defining skein relation for the Conway polynomial
implies that
$$
\nabla(\kbt) - \nabla(\kbo)
\ =\ -t\nabla(\kbs)
\ =\ \nabla(\kbo) - \nabla(\kbz),
$$
where the braids corresponding to the words $w'_1$ and $w'_2$ 
are depicted as rectangular boxes.

Therefore,
$$
\nabla(\kbt)-2\nabla(\kbo)+\nabla(\kbz)=0,
$$
as required.

(4).
We will prove that for any word $w$ we have $\chi(wBC)=t^2\chi(wB)$.
Indeed,
\begin{eqnarray*}
\chi(wB)&=&\sum_{w'\subseteq w}(-1)^{|w|-|w'|}
    \bigl(\nabla(w'B)-\nabla(w')\bigr),\\
\chi(wBC)&=&\sum_{w'\subseteq w}(-1)^{|w|-|w'|}
    \bigl(\nabla(w'BC)-\nabla(w'B)-\nabla(w'C)+\nabla(w')\bigr).
\end{eqnarray*}

Now, the needed assertion follows from the identity
$$
  \nabla(w'BC)-\nabla(w'C) = (t^2+1)(\nabla(w'B)-\nabla(w')),
$$
which is proved by applying the Conway skein relation:
\newcommand{\lfo}{\rb{-12mm}{\ig[height=27mm]{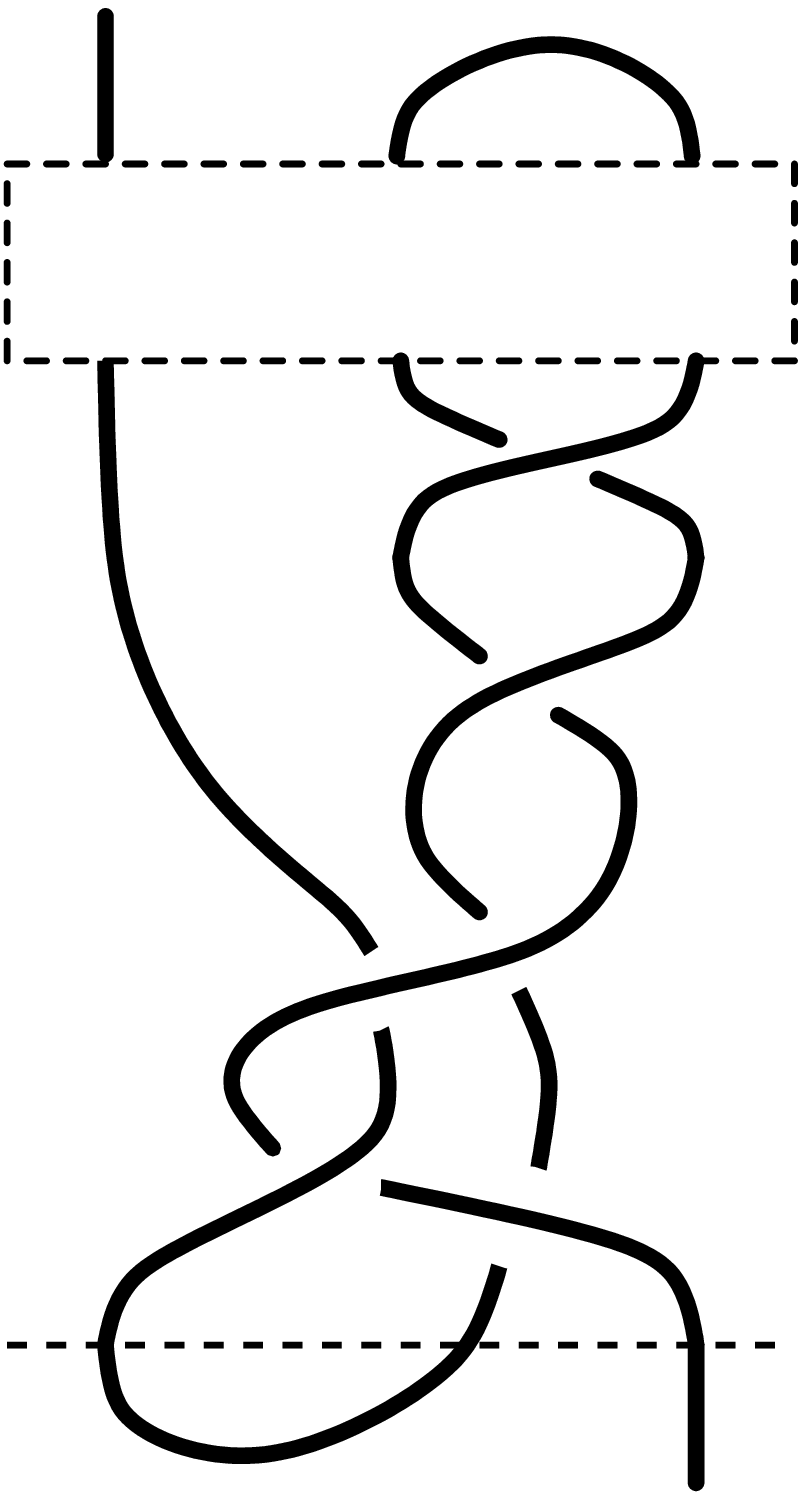}}}
\newcommand{\lfd}{\rb{-12mm}{\ig[height=27mm]{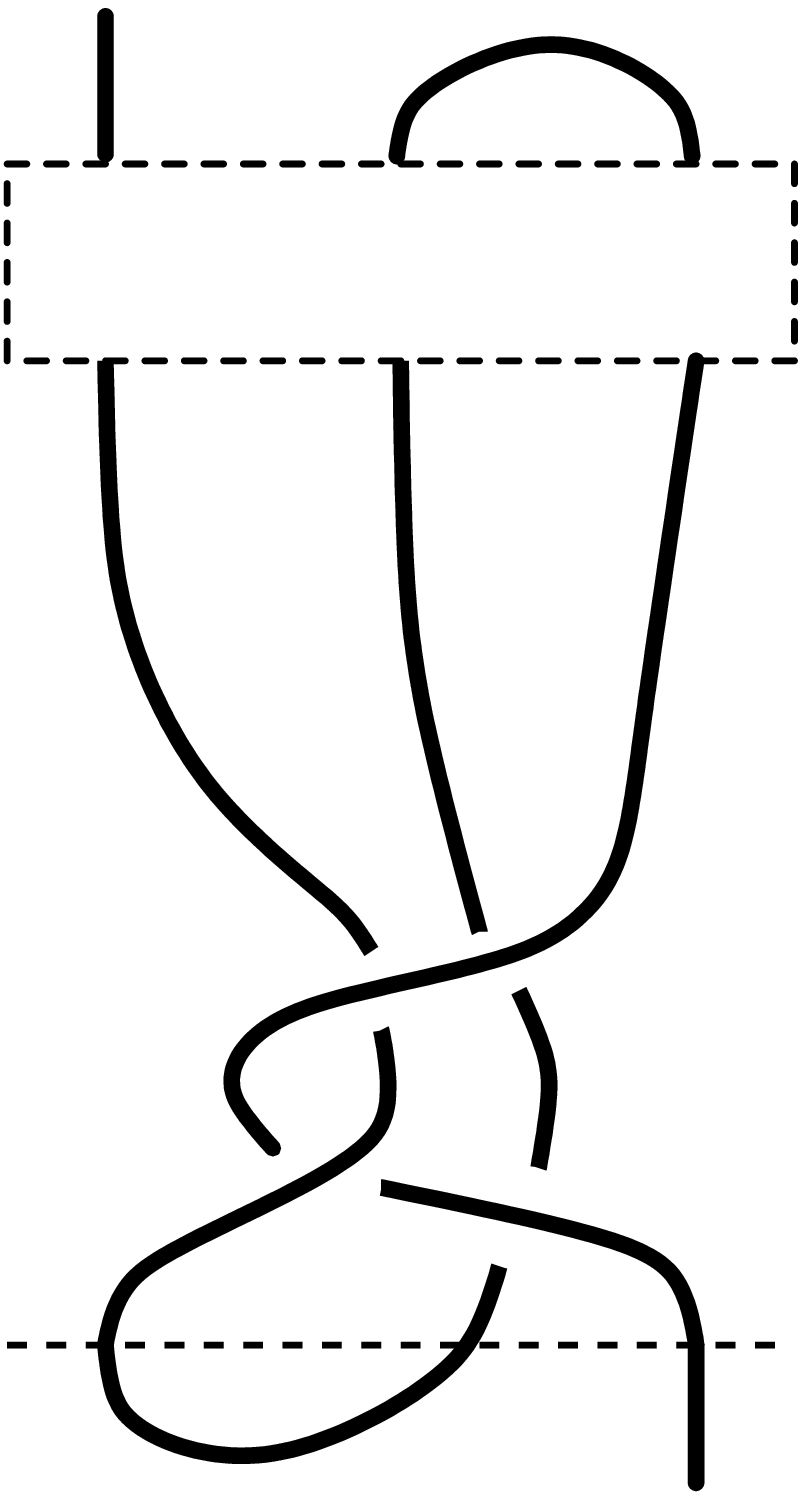}}}
\newcommand{\lft}{\rb{-12mm}{\ig[height=27mm]{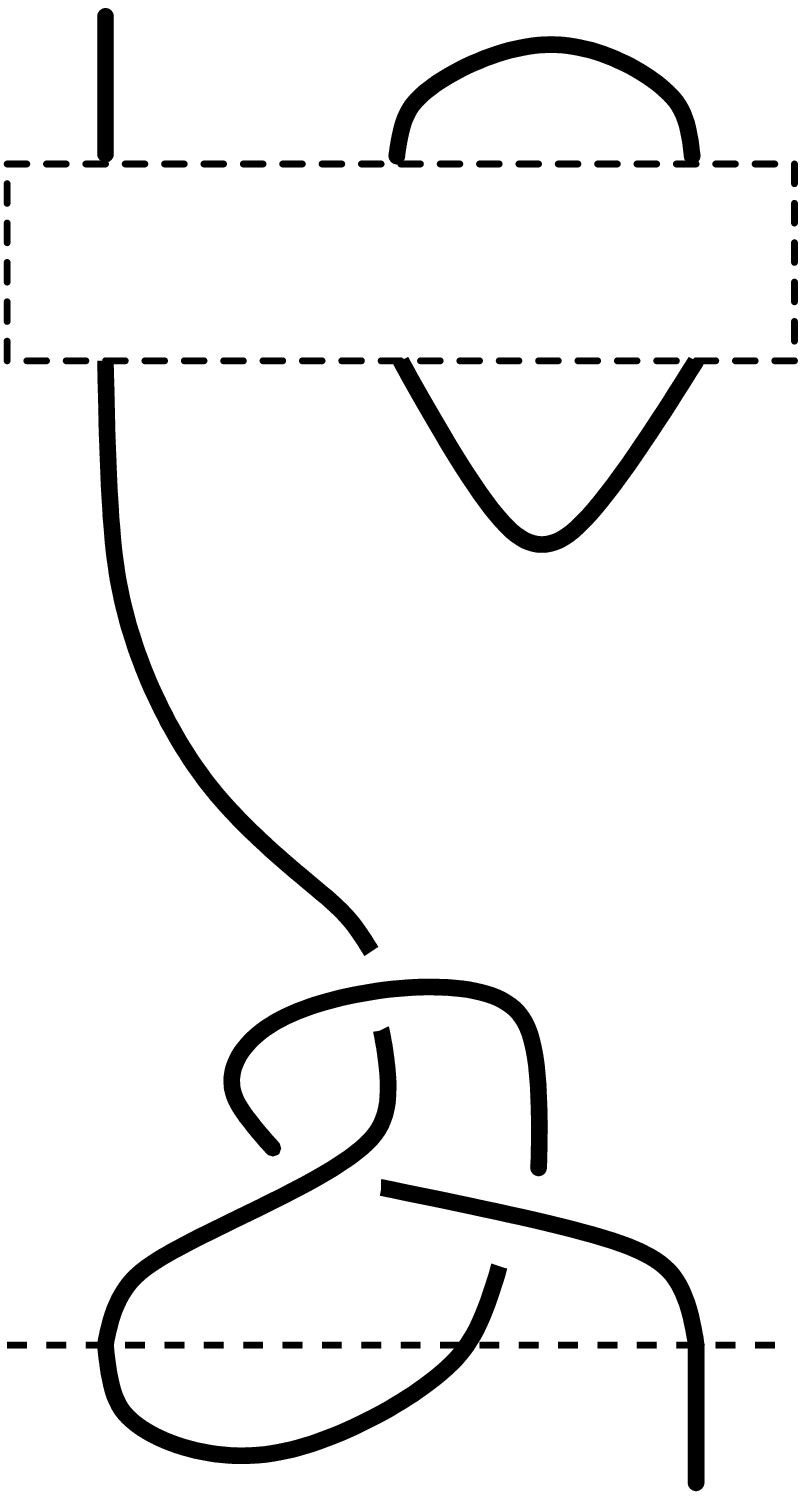}}}
\newcommand{\lfc}{\rb{-10mm}{\ig[height=23mm]{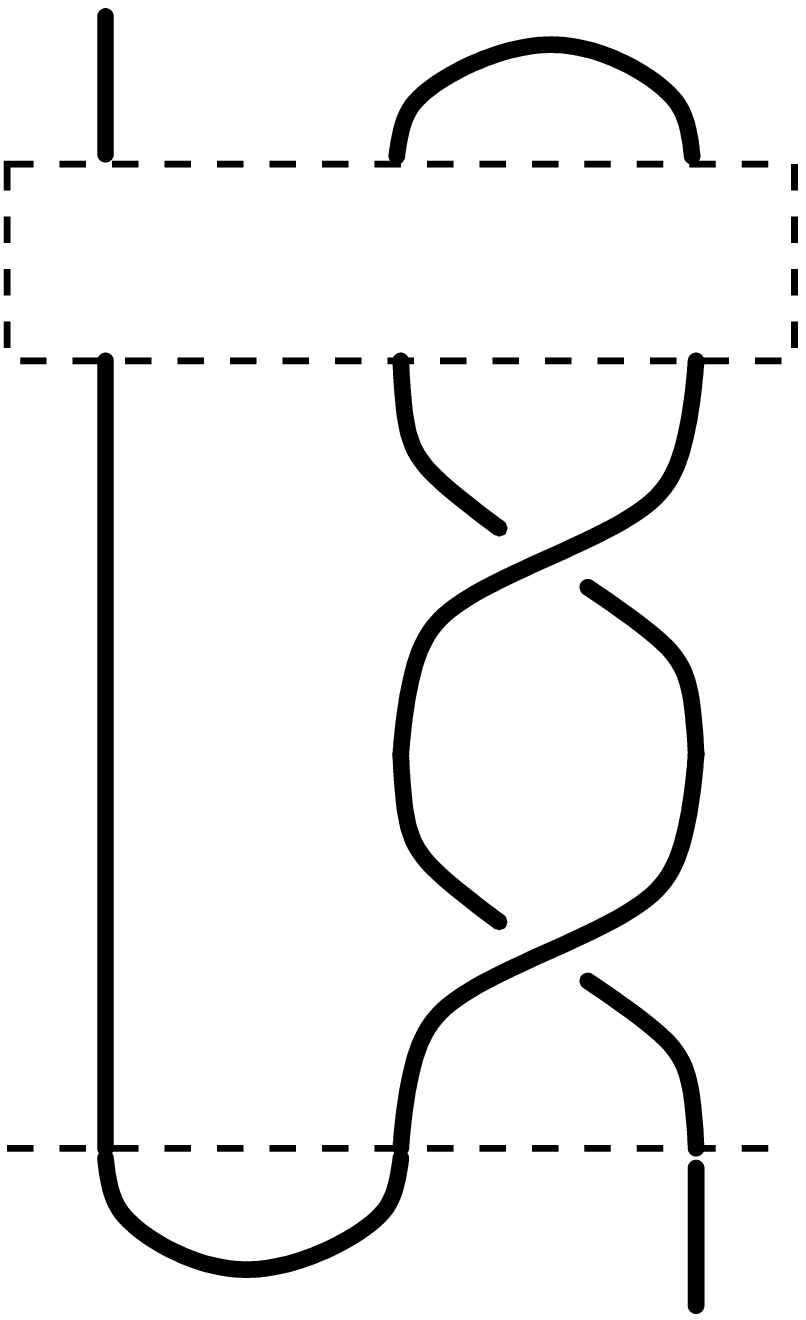}}}
\newcommand{\lfp}{\rb{-10mm}{\ig[height=23mm]{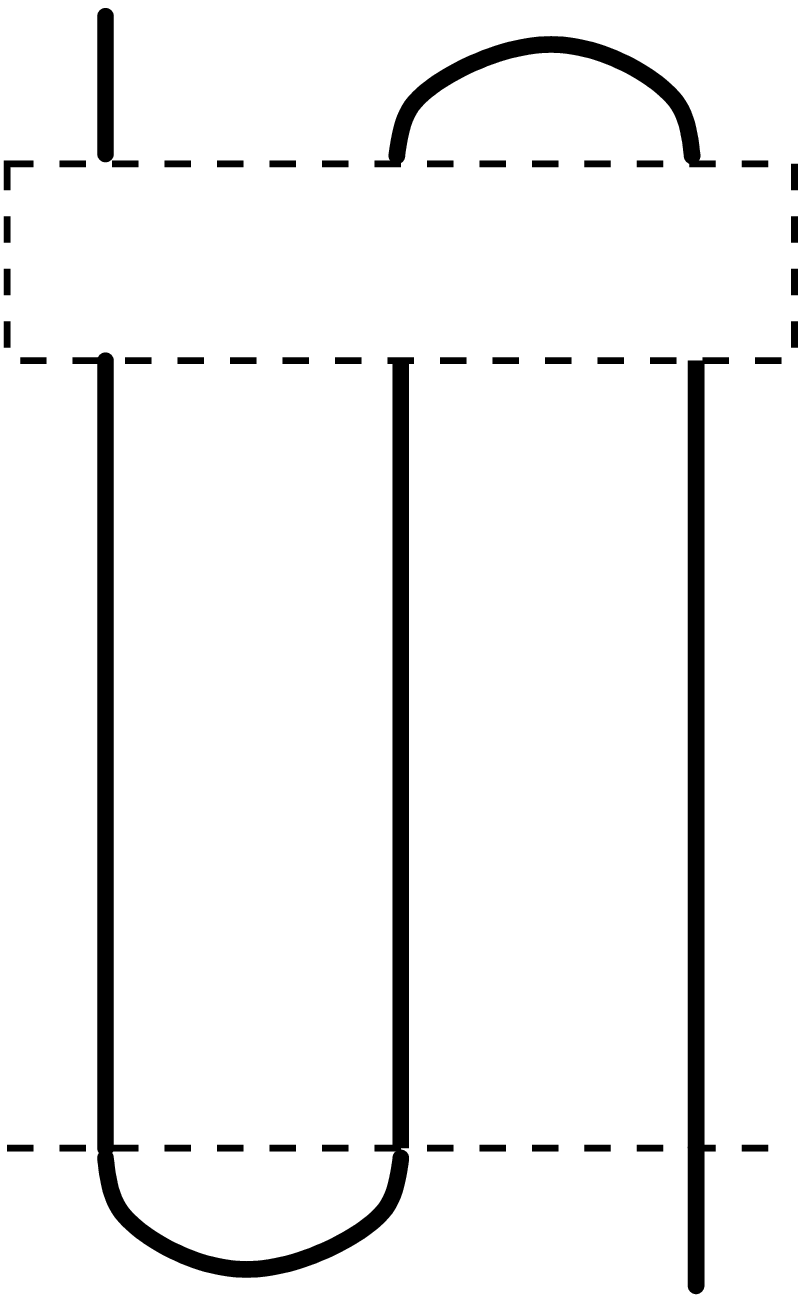}}}
\newcommand{\lfs}{\rb{-10mm}{\ig[height=23mm]{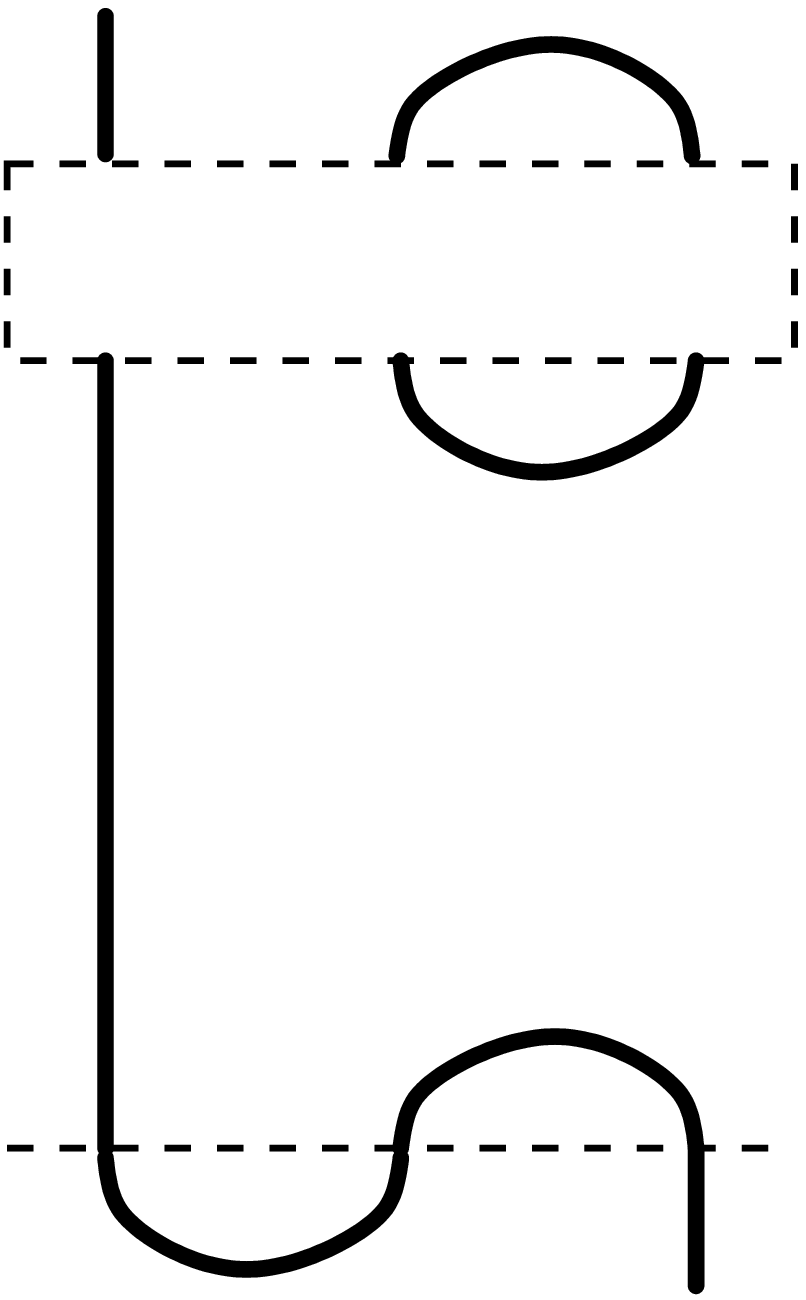}}}
$$
  \nabla(\lfo) - \nabla(\lfd) = -t \nabla(\lft),
$$
$$
  \nabla(\lfc) - \nabla(\lfp) = -t \nabla(\lfs)
$$
and noticing that taking the sum of a link with the trefoil knot
(which happens in the pictures on the right) leads to the multiplication of
the corresponding polynomial by $(t^2+1)$.
\smallskip

(5).
Here we must show that
$\chi(C^nBw) = -p_1p_{n-1}\chi(w)$
where $w$ is an arbitrary word in $C$ and $B$.
Indeed, let us split the alternating sum for $\chi(C^nBw)$ into the parts 
corresponding to a fixed subword $w'\subseteq w$:
$$
\chi(C^nBw)=\sum_{w'\subseteq w}
(-1)^{|w|-|w'|}\sum_{l=1}^n(-1)^{n-l}\binom{n}{l}
\bigl(\nabla\k(C^lBw')-\nabla\k(C^lw')\bigr).
$$
Using the Conway skein relations on a proper crossing, we get:
$$
  \nabla\k(C^lBw')-\nabla\k(C^lw')\, =\, -t\nabla(K_{(l)}),
$$ 
where $K=\k(w')$ and
$K_{(l)}$ denotes the oriented 2-component link obtained from the oriented
knot $K$ by adding a trivial $l$-linked component to $K$ according to the
picture (example for $l=3$):
$$
  \ig[height=22mm]{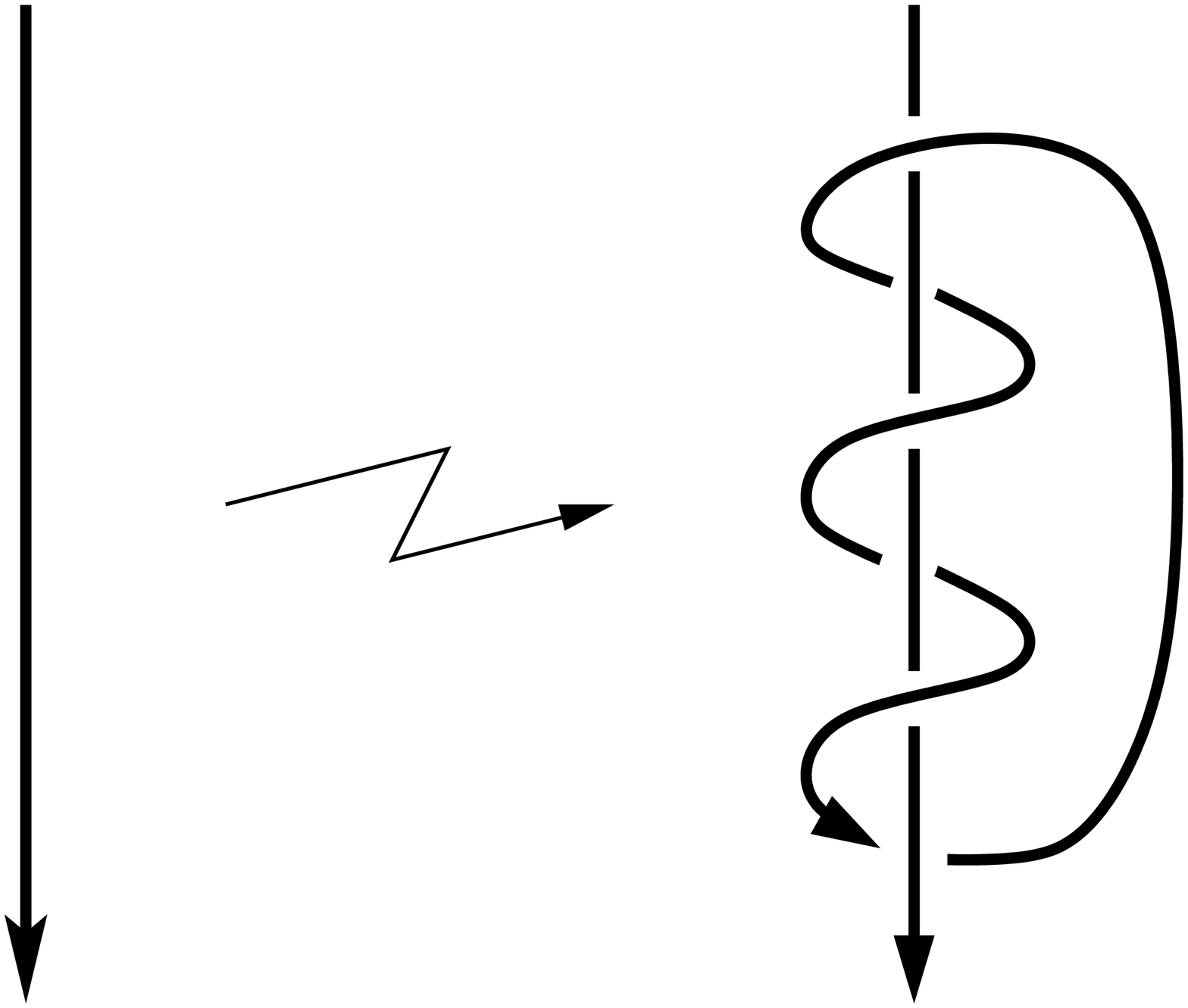}
$$

{\bf Lemma.} For any knot $K$ and any natural number $l$, we have
$$
  \nabla(K_{(l)})=t(q_0+q_1+...+q_{l-1})\nabla(K),
$$
where $q_s$ is the Conway polynomial of the torus knot of type $(2,2s+1)$
given explicitly by 
$$
  q_s=\sum_{j=0}^s\binom{s+j}{s-j}t^{2j}.
$$

This equality is proved, as usual, by recursively applying the Conway skein 
relation.
Substituting it into the previous formula for $\chi(C^nBw)$, we get
$$
  -t^2\cdot\sum_{l=1}^n(-1)^{n-l}\binom{n}{l}\sum_{s=0}^{l-1}q_s\cdot
   \sum_{w'\subseteq w}(-1)^{|w|-|w'|}\nabla(\k(w')),
$$
and it remains to show that the middle term of this product is equal to
$p_{n-1}$. Indeed, it is easily transformed to the form
$$
  \sum_{s=0}^{n-1}(-1)^{n-1-s}\binom{n-1}{s}q_s
$$
or, recalling the expression for $q_s$, to
$$
  \sum_{s=0}^{n-1}\sum_{j=0}^s(-1)^{n-1-s}
  \binom{n-1}{s}\binom{s+j}{s-j}t^{2j}.
$$
Changing the order of summation, this can be rewritten as
$$
  (-1)^{n-1}\sum_{j=0}^{n-1}\Bigl[\sum_{s=j}^{n-1}
  (-1)^s\binom{n-1}{s}\binom{s+j}{2j}\Bigr]t^{2j}.
$$

An application of the product summation formula
(5.24) from \cite{GKP} to the sum over $s$ inside the brackets
gives $(-1)^{n-1}\binom{j}{2j-n+1}$ thus proving the required assertion.

\smallskip

{\bf Remark.} The coefficients of the polynomials $p_n$ and $q_n$ can be 
read off the Pascal triangle in this way:
$$
  \ig[height=50mm]{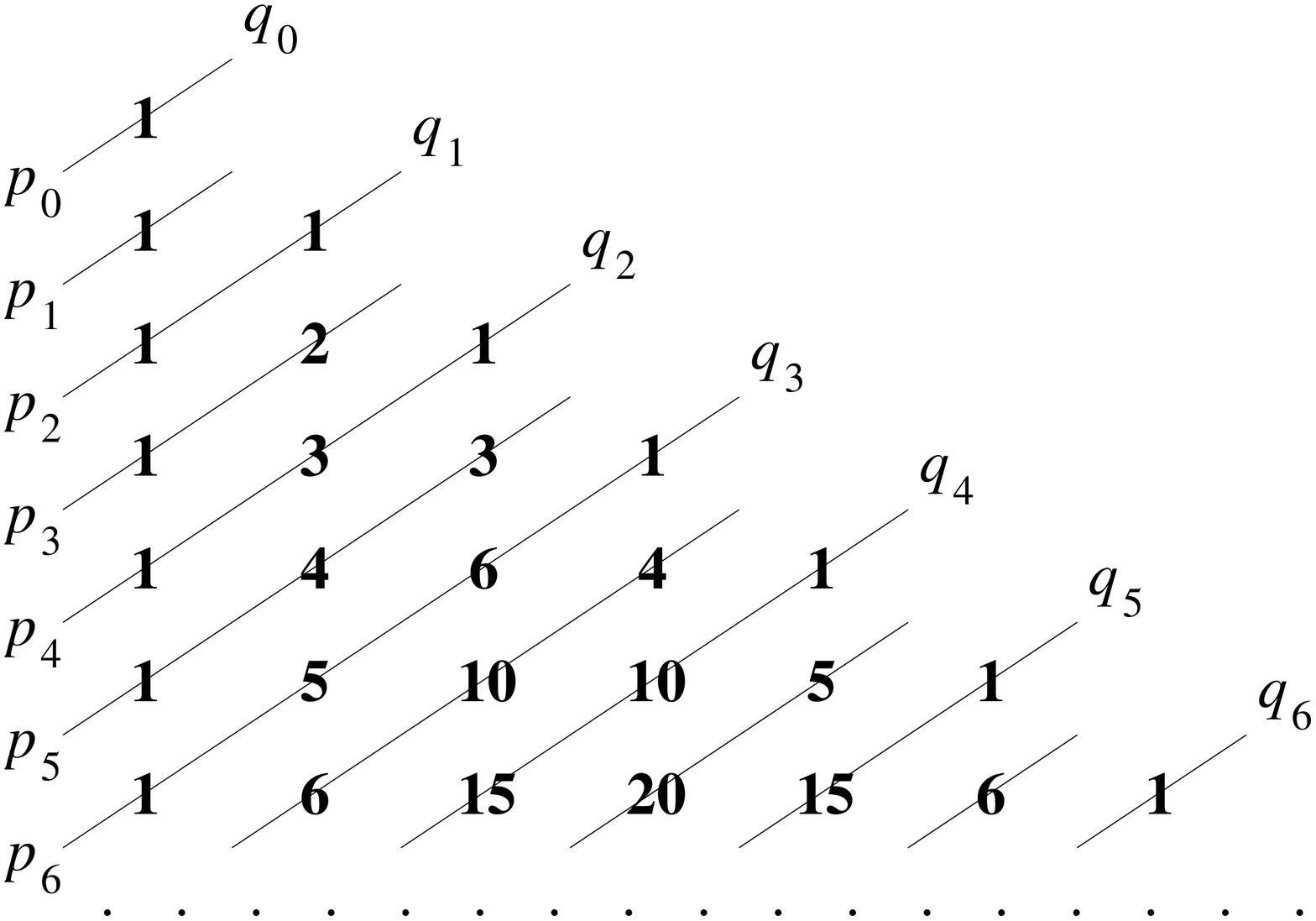}
$$

\section{Evaluation on the associator}
\label{sec-assoc}

The Drinfeld associator \cite{Dr,LM} is a remarkable 
element of the algebra $\At$, given by an infinite series 
in the (non-commuting) variables $a=A/(2\pi i)$, $b=B/(2\pi i)$ with 
coefficients in the algebra of multiple zeta values (MZV, see \cite{Hoff}):
\begin{eqnarray*}
 \Phi &=&1-\z_2[a,b]-\z_3([a,[a,b]]+[b,[a,b]])\\
     && -\,\z_4[a,[a,[a,b]]]-\z_{3,1}[b,[a,[a,b]]]-\z_{2,1,1}[b,[b,[a,b]]]+
         \frac{1}{2}\z_2^2[a,b]^2...
\end{eqnarray*}
(see \cite{Du} for an explicit expansion of $\Phi$ up to degree 12).

Taking the value of the complexified mapping 
$\chi_\C:\hat\A^h(3)\otimes\C\to\C[t]$ on this element
$\Phi$, we obtain the following result:

\begin{conjecture}
\begin{eqnarray*}
\chi_\C(\Phi) &=& -\z_2T^2 +(-\z_3+\z_{2,2})T^4
+(-\z_4+\z_{2,3}+\z_{3,2}-\z_{2,2,2})T^6+...\\
&=&\sum_{n=1}^\infty\big(\sum_{k=1}^n(-1)^k\z_{n+k}^{(k)}\bigr)T^{2n},
\end{eqnarray*}
where $T=t/(2\pi i)$, the numbers $\z_{l_1,...,l_k}=\z(l_1,...,l_k)$ 
are multiple zeta values
and $\z_m^{(k)}$ is the short-hand notation for the sum of
all $\z_{l_1,...,l_k}$ where $l_i\ge2$ for $i=1,...,k$
and $l_1+l_2+...+l_k=m$.
\end{conjecture}

We have checked this formula by computer up to $T^{10}$ (see \cite{Du}) using 
the table of relations between MZV's provided in \cite{Pet}.
It is an interesting remark that, when evaluated numerically, the
coefficients of this polynomial:
\begin{eqnarray*}
 &&-1.644934\,T^2 -0.390314\,T^4 -0.332698\,T^6 -0.312405\,T^8 \\
 &&-0.303958\,T^{10}-0.300153\,T^{12} -0.298365\,T^{14}-0.297505\,T^{16}+\dots
\end{eqnarray*}
seem to tend to a limit whose nature remains unclear.
 
\section{Open problems}
\label{sec-prob}

\begin{enumerate}
\item
For what triples of reduced rational fractions $p_1/q_1$, $p_2/q_2$, 
$p_3/q_3$ where both numerators and denominators form arithmetical 
progressions,
the values of the Conway polynomial on the corresponding rational knots
also form an arithmetical progression?
Our proof of assertion 3 of the main theorem (see page \pageref{mt3})
gives an abnormally big number of such triples, for instance, 
$(\frac{11}{3},\frac{31}{13},\frac{51}{23})$, 
$(\frac{11}{5},\frac{19}{9},\frac{27}{13})$, 
$(\frac{17}{5},\frac{41}{11},\frac{65}{17})$, 
$(\frac{75}{13},\frac{111}{19},\frac{147}{25})$, 
yet the claim about arithmetical progressions is not true in general.

\item
Generalize the main theorem (Section \ref{sec-scp}) in two directions:
(A) to pure braids with an arbitrary number of strands, (B) to the HOMFLY
polynomial which is a generalization of the Conway polynomial.

\item
This is related to the Remark on page \pageref{magnus_gen}.
Describe all (or some) triples of formal series $P,Q,R$ in $\At\otimes\C$
starting with terms of degree higher than one,
such that the correspondence $x_{12}\mapsto 1+t_{12}+P$, $x_{23}\mapsto
1+t_{23}+Q$, $x_{13}\mapsto 1+t_{13}+R$ defines a group homomorphism. 

\item
Prove the conjecture from Section \ref{sec-assoc}.
Find similar facts for other horizontal associators and relate them
to the action of the Grothendieck--Teichm\"{u}ller group (see \cite{Dr,ES}).

\end{enumerate}


\begin{thebibliography}{JM}

\bibitem{BN} D.\,Bar-Natan, {\it Vassiliev and Quantum Invariants 
of Braids}, Proceedings of Symposia in Applied Mathematics {\bf 51} 
(1996) 129--144. Online at {\tt arXiv:q-alg/9607001}.

\bibitem{Bir} J.\,S.\,Birman, {\it Braids, Links and Mapping Class
Groups}, Princeton University Press, 1974.

\bibitem{CDM} S.\,Chmutov, S.\,Duzhin and J.\,Mostovoy. {\it CDBooK. 
Introduction to Vassiliev Knot invariants}, draft version of a book,\\
Online at \verb#http://www.pdmi.ras.ru/~duzhin/papers/cdbook/#.

\bibitem{Dr} V.\,Drinfeld, {\it On quasi-triangular quasi-Hopf
    algebras and a group closely connected with
    $\mathrm{Gal}(\overline{\Q}/\Q)$}, Algebra i Analiz,
{\bf 2}, no.~4, 149--181 (1990). English translation:
Leningrad Math. J., {\bf 2} (1991)  829--860.

\bibitem{Du} S.\,Duzhin, {\it Program and data files related to the Drinfeld
associator}, online at
\verb#http://www.pdmi.ras.ru/~arnsem/dataprog/associator/#.

\bibitem{ES} P.\,Etingof, O.\,Schiffmann,
{\it Lectures on Quantum Groups}, International Press, 1998.

\bibitem{GKP} R.\,Graham, D.\,Knuth and O.\,Patashnik. 
{\it Concrete Mathematics},
Reading, Massa\-chusetts: Addison-Wesley, 1994.

\bibitem{Hoff}
M.\,E.\,Hoffman, {\it Multiple harmonic series}, Pacific J. Math. 152 (1992),
275-290. 

\bibitem{LM} T.\,Q.\,T.\,Le, J.\,Murakami, {\it The Kontsevich integral
   for the Kauffman polynomial},
   Nagoya Mathematical Journal {\bf 142} (1996) 39--65.

\bibitem{MSt} J.\,Mostovoy, T.\,Stanford,
    {\it On a map from pure braids to knots},
 Journal of Knot Theory and its Ramifications, {\bf 12} (2003), 417--425.\\
Online at \verb#http://www.matcuer.unam.mx/~jacob/works.html#.

\bibitem{MW} J.\,Mostovoy, S.\,Willerton,
    {\it Free groups and finite-type invariants of pure braids.}
    Math.\ Proc.\ Camb.\ Philos.\ Soc.\ {\bf 132} (2002) 117--130.\\
Online at \verb#http://www.matcuer.unam.mx/~jacob/works.html#.

\bibitem{Mur} K.\,Murasugi, {\it Knot Theory and Its Applications},
Birkh\"{a}user, 1996.

\bibitem{Pap} S.\,Papadima,
{\it The universal finite-type invariant for braids, with integer
coefficients,} Topology Appl.\ {\bf 118} (2002) 169--185.

\bibitem{Pet} M.\,Petitot, {\it Tables of relations between MZV up to weight
16},\\ Online at \verb#http://www2.lifl.fr/~petitot/#.

\bibitem{PS} V.\,V.\,Prasolov, A.\,B.\,Sossinsky, 
   {\it Knots, Links, Braids and 3-Manifolds}. 
   Translations of Mathematical Monographs, {\bf 154}.
   American Mathematical Society, Providence (1997).

\end{thebibliography}
\end{document}